\documentclass[a4paper,11pt,reqno]{amsart}

\usepackage{amsmath,amssymb}
\usepackage{enumerate}
\usepackage{ifthen}
\usepackage{calc}
\nonstopmode\numberwithin{equation}{section}
\newtheorem{definition}{Definition}[section]
\newtheorem{theorem}{Theorem}[section]

\newtheorem{remark}{Remark}[section]

\allowdisplaybreaks

\begin{document}
\title{ group ring based  public key cryptosystems}
\maketitle
 \begin{center}
\center{  $\text{Gaurav Mittal}^{1, 2}, \text{Sunil Kumar}^2,\text{Shiv Narain}^3$, and $\text{Sandeep Kumar}^2$ }
\medskip
{\footnotesize

 \center{$^1$Department of Mathematics, Indian Institute of Technology Roorkee, Roorkee, India\\$^{2}$DRDO, India \\$^{3}$Department of Mathematics, Arya P. G. College, Panipat, India\\ emails: gmittal@ma.iitr.ac.in, sunil\_kumar@hqr.drdo.in,\\  math.shiv@gmail.com,  sandeepkumar@hqr.drdo.in

 }
 }

\bigskip
\begin{abstract} In this paper, we propose two cryptosystems based on group rings and  existing cryptosystem. First  one is Elliptic ElGamal type group ring public key cryptosystem whose security is greater than security  of cryptosystems based on elliptic curves discrete logarithmic problem (ECDLP). Second  is ElGamal type  group ring public key cryptosystem, which is analogous to ElGamal public key cryptosystem but has comparatively greater security. Examples are also given for  both  the proposed cryptosystems.
\end{abstract}
\end{center}
 \subjclass \textbf{Mathematics Subject Classification (2010)}: {94A60, 20C05, 20C07}

\keywords{\textbf{Keywords:}  group ring, units, cryptography, public key cryptography
\section{Introduction}

 Cryptography is the  study of mathematical techniques for  achieving Privacy, Authentication, Integrity  and
Non-repudiation, i.e. PAIN. For the secrecy and authenticity of sensitive information, cryptosystems are used extensively.  The one and only aim of these cryptosystems is to ensure end-to-end communication between two or more parties such that the information encrypted at one end is understood only at the other end. For  end-to-end secure communication, various symmetric ciphers have been developed till $1970$'s, for instance, classical substitution ciphers, transposition ciphers, rotor machines (see $[20]$) etc.   Since then many symmetric ciphers came in the picture but security-wise some of the best known     are   DES  and AES both based on the concepts of confusion and diffusion introduced by Claude Shannon $[17]$,  however, in the current scenario DES is no longer in use because the key length does not withstand against the brute force attack.

Some  of the greatest revolutions in  cryptography are  the development of asymmetric or public key cryptography and key authentication protocols. 
The well known public key cryptosystesms or ciphers till date are   RSA, ElGamal, NTRU  (see $[9, 20]$) etc. Keeping the keys out of the reach of foe is extraordinarily imperative in public key cryptosystems. This leads to the concept of 
key authentication protocols. 
 Various key authentication protocols have been developed by many researchers in the past (cf. [20] for more details). Very recently, Chandrashekhar et al. [6] proposed a new efficient and secure  authentication protocol for discrete logarithm define over multiplicative group based cryptosystem using the basic methodology.

  The security of most of the  public key cryptosystems as well as key authentication protocols depend on the hardness of solving some  underlying mathematical problem with a trapdoor like integer factorization, discrete logarithmic, shortest non-zero vector in the lattice etc. and most of these problems arise in number theory $[9]$. 
 Combinatorial group theory has also attracted much interest for the construction of public key cryptosystems $[1, 8, 19]$. The hard underlying mathematical problem associated with most of these cryptosystems is   conjugacy search problem which can be considered as the generalization of discrete logarithm problem to groups other than $\mathbb{F}_p$, where $p$ is some prime.
 
Hill cipher and Affine Hill cipher  [20] are well known classical ciphers extensively used  in symmetric cryptography. But, now a days these ciphers are not only limited to symmetric key cryptography. 
 Using the classical Affine  Hill cipher, a public key cryptosystem has been developed in    [21] based on modulation of prime number. We also refer to references in [21] for more literature on the public key cryptosystems that are developed using symmetric cryptosystems. 
 
  Elliptic curves over finite fields  are very important from the cryptographic point of view (cf. $[12, 13]$). The interesting problem of adding two points on an elliptic curve was a very subtle issue and theory behind it leads to the fascinating field of mathematics known as Algebraic Geometry. For additional reading devoted to elliptic curve cryptography, see $[4, 14]$.

   To the best of our knowledge, algebraic structure known as group rings has not been used in cryptography upto first decade of $21^{\text{st}}$ century. Hurley $[11]$ discussed the various applications of group rings to the areas of communication.  In $[10]$,  the first known symmetric key cryptosystem using group rings was demonstrated with  few examples which show  how to combine RSA and unit of group ring and how to use  large power of a unit  as a public key   but no public key cryptosystem was discussed explicitly.

In this paper, we discuss two novel group rings based public key cryptosystems. First one is  Elliptic  Elgamal type group ring public key cryptosystem and second one is ElGamal type group ring public key cryptosystem  without the involvement of Elliptic curves. 
Elliptic ElGamal public key cryptosystem which is analogous to  ElGamal is discussed in $[9]$. The main difficulty with this cryptosystem is non-obvious way of attaching the  plaintext messages and elements of  elliptic curve over finite field. We try to resolve this problem in our work by  combining the elements of group ring with that of elliptic curves in a unique manner. The sounding effect of these proposed public key cryptosystems is their enhanced security  in comparison to that of existing cryptosystems. This enhanced security can be well understood  for traditional classical computers but in the coming era of quantum computers where public key cryptosystems (in particular, RSA, ECDSA) are no longer safe $[3, 16, 18]$, this can play a vital role because of the involvement of two operations in group ring which are connected via distributive law. 

Rest of the paper  is designed as follows:  Section $2$ pertains to  preliminary definitions and discussion of both the novel cryptosystems. Examples related to these schemes are discussed in Section $3$. Section $4$ is related to discussion of availability of units of group rings and last section comprises  some concluding remarks.

\section{Cryptographic systems in group rings}
In the present section, we propose two public key cryptosystems in group rings. But before that, we provide some preliminary definitions and results  requisite for our work.  
\begin{definition}
Group Ring: Let $R$ be a ring with unity and $G$ be a group. The set $$RG= \bigg\{\sum_{j=1}^t r_jg_j: r_j\in R, g_j \in G\bigg\}$$of all finite sums is known as  group ring. As name suggests, RG is a group, ring, and module over $R$. See $[15]$ for more information on group rings. For $a, b \in RG$, $a*b$ denotes the product of $a$ and $b$.
\end{definition}
\begin{definition}
Unit of  a group ring: Let $u\in RG$. This $u$ is a unit of $RG$ if there exists an element $v\in RG$ such that $$u*v=v*u=1.$$Set of all units of $RG$ forms a group under the operation $*$.
\end{definition}
\begin{definition}
Normalized unit of  a group ring: Let $u= \sum_{i=1}^k\alpha_ig_i\in RG$ with $\alpha_i \in R$ and $g_i \in G$ be a unit of $RG$. This $u$ is said to be normalized if $$\alpha_1+\alpha_2+\cdots+\alpha_k=1.$$Set denoted by $V(RG)$ of all the normalized units of $RG$ forms a group under the operation $*$.
\end{definition}
\begin{definition}
 Discrete Logarihtmic Problem $($DLP$)$: 
 Let $g, h$ be two  elements of $\mathbb{F}_p$ such that $h$ is some power of $g$. Then DLP is the problem of finding an integer $n$ such that $h=g^n$.
\end{definition}
\begin{definition}
Discrete Logarihtmic Problem $($DLP$)$ in group ring: Let $u$ be an element of $RG$ and $v$ be another element of $RG$ which is some power  of $u$. Then, DLP in group ring is determination of an integer $x$ such that $v=u^x.$ 
\end{definition}
\begin{definition}
Elliptic Curve: Consider the Weierstrass equation $$E:Y^2= X^3+RX+S.$$ Set of all the solutions of $E$ together with an extra point $O$ that lives at infinity is an elliptic curve. Here $R$ and $S$ are constants which satisfy $4R^3+27S^2\neq 0.$ 
\newline 
Further, if we consider only those elements of $E$ which also belong to   $\mathbb{F}_p\times \mathbb{F}_p$, then the set
$$E(\mathbb{F}_p)=\big\{(x, y): (x, y) \in \mathbb{F}_p\times \mathbb{F}_p\ \text{and}\ y^2=x^3+Rx+S\big\}\cup\{O\}$$
is an elliptic curve over the finite field $\mathbb{F}_p$ with $R, S \in  \mathbb{F}_p$.
\end{definition}
\begin{definition}
Elliptic Curve Discrete Logarihtmic Problem $($ECDLP$)$: Let $P, Q$ be two  points on elliptic curve $E(\mathbb{F}_p)$. Then ECDLP is the problem of determination of an integer $n$ such that $Q=nP$,  provided such an integer exists.
\end{definition}
Next, we give the explicit formulas to add and subtract any two points on an elliptic curve. For proof see $[9, \text{Theorem} \ 5.6]$.
\begin{theorem}
Let $P$ and $Q$ be any two points on the elliptic curve $E(\mathbb{F}_p)$. \begin{enumerate}
\item If $P=O$, then $P+Q=Q$. Similarly, if $Q=O$, then $P+Q=P$. \item Otherwise, write $P=(x_1,y_1)$ and $Q=(x_2,y_2)$.\item If $x_1=x_2$ and $y_1+y_2=0$, then $P+Q=O$.\item Otherwise, define $\lambda$ by $$\lambda=\begin{cases} \dfrac{3x_1^2+R}{2y_1}\quad \text{if}\ P=Q\\ \dfrac{y_2-y_1}{x_2-x_1}\quad \text{if}\ P\neq Q.\end{cases}$$ Then $P+Q=(x_3, y_3)$ with $$x_3=\lambda^2-x_1-x_2\qquad\text{and}\qquad y_3=\lambda(x_1-x_3)-y_1.$$
\end{enumerate}

\end{theorem}
Now, we discuss the important task of associating each data or number to be encrypted with the elements of group ring $RG$ where $R$ is the ring with identity and  $G= \{g_1, g_2, \cdots\}$ is either finite or countable group. In our study, we either consider $R=\mathbb{F}_p$ for some  prime $p>0$ or $R=\mathbb{Z}$.
\subsection{Connection between given data  and  elements of group ring}
 Suppose that all the digits of the  data to be encrypted  are elements of $R$, for instance, data consists of elements ranging from $0$ to $p-1$ for some prime $p$ and $R=\mathbb{Z}_p$.   Given any  arbitrary data, write it in blocks of length $t$ where $t\in \mathbb{Z}^+$. Further, any arbitrary   data block  $m=m_1m_2\cdots m_t$ with $m_i\in R$  can be considered  as an element  of group ring $RG$ via the following representation $$r=\sum_{i=1}^tm_ig_i, \ \ g_i \in G.$$  For above representation to be unique, order of $G$ must be atleast $t$ and therefore we assume the same, i.e. $|G|\geq t$ and can be atmost countable.  
 If  length of the message is not a multiple of $t$, then padding can be done with $0'$s.
\subsection{Elliptic  ElGamal type group ring public key  cryptosystem} Given a message $r=\sum_{i=1}^tm_ig_i$,  write it uniquely as a row vector $$r=\big[
m_1\ \ \ m_2\ \ \ \cdots \ \ \ \cdots \ \ \ m_t\big]
$$ where $m_i$ is the coefficient of element $g_i$ in $r$. 
Choose a point $P= (x_1, y_1)$ on the elliptic curve $E(\mathbb{F}_p)$  and a unit $u$ of  group ring $RG$ such that order of $u$ is  large. This unit is not made public. 
Further,  secretly choose two positive integers $(n_1, n_2)$ and use these integers to compute $$Q=n_1P= (x_2, y_2)\ \ \text{and}\ \ A=u^{n_2}$$ where $A\in RG$. 
Set $(P, Q, A)$ as the public key for encryption and $(n_1, u^{-n_2})$ as the private key for decryption.  Now  select a random one time key $n_3$ and let $n_3Q=(x_3, y_3)$. This key is used for encrypting only one message or some fixed number of blocks that needs to be decided in the beginning and then discarded.
\subsubsection{Encryption}
For encryption, message $r$ is encrypted using public key to obtain the ciphertext  $(C_1, C_2)$ where 
$$C_1=n_3P, \ \ C_2= (r\oplus n_3Q)*A$$ where the $\oplus$ operation is defined as $$r\oplus n_3Q= \big[
m_1\ \ \ m_2\ \ \ \cdots\ \ \ \cdots\ \ \ m_t
\big]\oplus (x_3, y_3)\hspace{45mm}$$  $$= \big[
m_1+x_3+y_3\ \ \ m_2+x_4+y_4\ \ \ \cdots\ \ \  \cdots\ \ \ m_t+x_{t+2}+y_{t+2}\big]
$$ with $$n_3Q=(x_3, y_3),\ 2n_3Q=(x_4, y_4), \ 3n_3Q=(x_5, y_5),\ \cdots, \ tn_3Q=(x_{t+2}, y_{t+2})$$ and since $m_i, x_i, y_i\in \mathbb{F}_p$, above addition is  defined. Rationale behind the $\oplus$ operation is to  combine an element of group ring with an element of elliptic curve such that the result is an element of group ring.  So, the ciphertext $(C_1, C_2)$ has two parts where $C_1$ is the element of elliptic curve and $C_2$ is the element of group ring. This ciphertext is then send to the receiver say Alice.
\begin{remark}
If the message $m$ is of the form $$\big[
m_1\ \ \ 0 \ \ \ m_3\ \ \  m_4\ \ \ 0\ \ \ \cdots\ \ \ 0\big],$$ then we can simply encrypt it in the form $\big[
m_1\ \ \ 0 \ \ \ m_3\ \ \ m_4\big]
$. But the $0$ entries in between the non-zero entries of the message  cannot be skipped and needs to be considered in order to  uniquely obtain the  plaintext from ciphertext.
\end{remark}
\begin{remark}
Suppose that the operation $\oplus$ is of the form $$r\oplus n_3Q= \big[
m_1\ \ \ m_2\ \ \ \cdots\ \ \ \cdots\ \ \ m_t\big]
\oplus (x_3, y_3)\hspace{35mm}$$  $$= \big[
m_1+x_3+y_3\ \ \ m_2+x_3+y_3\ \ \ \cdots\ \ \ \cdots\ \ \ m_t+x_{3}+y_{3}\big]$$ where we are adding $x_3+y_3$ in every message digit. This added block can leak some information about the message if there are consecutive zeros in the message. To avoid such a situation we added the multiples of $(x_3, y_3)$ in the message digits.
\end{remark}
\subsubsection{Decryption}
Since Alice knows  $u^{-n_2}$, she computes   $$C_2*u^{-n_2}= (r\oplus n_3Q)*A* u^{-n_2}= r\oplus n_3Q.$$  Alice also knows $n_1$ and she use this to obtain $$n_1C_1= n_1n_3P= n_3n_1P= n_3Q.$$ Thereafter, she computes the multiples of $n_3Q$  and perform the decryption operation $\ominus$,   i.e. $$(C_2*u^{-n_2})\ominus n_1C_1= (r\oplus n_3Q)\ominus (x_3, y_3)\hspace{60mm}$$ $$=\big[
m_1+x_3+y_3\ \ \ m_2+x_4+y_4\ \ \ \cdots\ \ \ \cdots\ \ \ m_t+x_{t+2}+y_{t+2}
\big]\ominus (x_3, y_3) $$ $$=\big[
m_1\ \ \ m_2\ \ \ \cdots\ \ \ \cdots\ \ \ m_t\big]
= \sum_{i=1}^tm_ig_i  =r$$ where $\ominus$ operation means subtracting $x_3+ y_3, x_4+ y_4, \cdots, x_{t+2}+y_{t+2}$ from the corresponding elements of  row vector.
\begin{remark}
We can also consider $R=\mathbb{Z}$ in Elliptic ElGamal type group ring public key cryptosystem in place of $\mathbb{F}_p$. In that case, whole encryption scheme remains same with the only difference that operation $\oplus$ involves usual addition instead of addition modulo $p$.
\end{remark}
\begin{remark}
For Elliptic ElGamal type group ring public key cryptosystem, we can also consider fields of different characteristics. For instance, we can consider $E(\mathbb{F}_p)$ and $R=\mathbb{Z}_q$ where $p$ and $q$ are some primes. If $p<q$, then the operation $\oplus$ involves  addition modulo $q$ and we directly add the elements of elliptic curves to that of group ring. If $p>q$, then first we need to convert the elements of elliptic curve to elements modulo $q$ and then apply  operation $\oplus$ which again involves  addition modulo $q$. An example related to this situation is discussed in Section $3$.
\end{remark}
\subsubsection{Logic behind the name}
Proposed encryption scheme is a novel approach involving the important  operation  of combining the elements of a group ring with those of Elliptic curve. The parameters involved are similar to ElGamal  and hence the cryptosystem is named as Elliptic ElGamal type group ring public key  cryptosystem. 
\subsubsection{Security analysis of the scheme} 
To obtain the plaintext from  ciphertext, adversary in between   needs to:
\begin{enumerate}
\item   solve ECDLP.
\item  determine inverse of power a unit  or unit  from its power. 
\end{enumerate}
In $\mathbb{F}_p^*$, there are algorithms for solving DLP and the fastest known algorithm  has subexponential running time  known as index  calculus method $[9]$. But there are no subexponential running time algorithms for solving ECDLP which makes elliptic curves very useful in cryptography. The fast known algorithm takes approximately $\sqrt{p}$ times to solve ECDLP (for example Pollard's $\rho$ method $[9]$), i.e. to say that there is no polynomial time algorithm to solve ECDLP.

Now we discuss the hardness of inverse computation problem in group rings.  There exist units which can be constructed by some formula (let's say known type units), for instance, unipotent units, central units, Bass cyclic units, Bicyclic units $[15]$. On taking power of a known type unit, resultant unit  may not be of known type and it becomes computationally infeasible to determine the inverse of a such a unit for large groups, i.e. to say that determination of unit from its power  introduces the difficulty of DLP in group rings (how much difficult DLP is in group rings is discussed in  subsection $2.3.3$). 
 Moreover, we can also combine the units of known type to get units of unknown type, for instance $A= (u_1u_2)^{n_1}$ where $u_1$ and $u_2$ are units of known type and it becomes computationally infeasible to obtain the inverse of such a unit without knowing the inverse of constituent units. So, in general we can say that it is safe to use power of units, power of product of units as the public keys.
 
  Hence, the novel scheme is comparatively more secure than the schemes whose security relies on solving ECDLP, e.g. Elliptic Elgamal public key cryptosystem.  In the current scenario, recommended key sizes (in bits) are $2048, 3072, 7680, \cdots$ for RSA and Diffie-Hellman whereas $224, 256, 384,\cdots$ for  the elliptic curves. This means elliptic curves have an edge over RSA and Diffie-Hellman in terms of key size despite the highly structured nature of $E(\mathbb{F}_p)$. If we consider an elliptic curve over a prime $p$ having size $224$ bits, then we can say that Elliptic ElGamal type group ring public key cryptosystem has a security equivalent to security of cryptosystems on $E(\mathbb{F}_{p'})$  where size of $p'$ is $224+k$ bits for some positive integer $k$ where $k$ depends on the  structure of the unit group of group ring.
 
 Now we introduce our second cryptosystem. For this cryptosystem too, parameters involved are similar to ElGamal (although  cryptosystem is entirely different) and hence we name it as ElGamal type group ring public key  cryptosystem. 
\subsection{ElGamal type group ring public key  cryptosystem} 
Secretly choose a unit $u\in RG$  of large order and two positive integers $n_1$ and $n_2$. Then find $A_1=u^{n_1}$. Choose another unit $v\in RG$ of large order (open in  public domain) and set $A_2=v^{n_2}$. Declare $(A_1, A_2, v)$ as the public key and $(A_1^{-1}, n_2)$ as the corresponding private key.  Now as in ElGamal, select a random ephemeral key $k$. \subsubsection{Encryption}
For encryption, find the ciphertext  $(C_1, C_2)$ where 
$$C_1=v^k, \ \ C_2= (r*A_1)*A_2^k.$$ This ciphertext is then send to  Alice. 
\subsubsection{Decryption}
As the second part $n_2$ of private key is known to Alice, she finds $(C_1)^{-n_2}$ and use this to obtain  $$C_2*v^{-n_2k}= r*A_1.$$On multiplying above with $A_1^{-1}$ (first part of private key), Alice gets the message $r$.
\subsubsection{Security of the scheme} To obtain the plaintext from  ciphertext, Eve needs  $n_2$ and $A_1^{-1}$ which requires the determination of unit from its power and inverse of a given unit respectively. In other words,  Eve's first task is to solve the DLP in group ring and second task is inverse computation problem which is already discussed in  Subsection $2.2.4$. So,  we only  discuss the hardness of DLP in group rings. For that, we need to measure the number of necessary operations required to efficiently compute the discrete log. Let $G$ be an arbitrary group and $g\in G$ with $|g|=n$ for   $n \in \mathbb{Z}^+$. If computing $g^x$ for some $x \in \mathbb{Z}^+$ is treated as a single operation, then discrete log can be computed in at most $n$ operations via brute force attack. If $n$ is a $r$ $(\geq 80)$ bit number, then on an average, brute force attack takes $O(2^r)$ operations to compute discrete log which is exponential time and therefore is no longer practical. Now on considering the same situation of discrete log computation in group rings, one can easily see that brute force is again impractical  because a group ring involves two operations connected via distributive law whereas a group involves only one operation.  Further, we know that index calculus method is the  fastest known subexponential algorithm to solve DLP in groups (already discussed in Subsection $2.2.4$), but no such algorithm is known to solve DLP in group rings. However, one thing we can  claim that currently there exists no  polynomial time algorithm to solve DLP in group rings because if there is some algorithm, then it is also implementable on groups. But the best known for groups is index calculus method.

So, we can say that currently there is no threat to DLP from classical computers for properly chosen groups as well as group rings. But 
 there exist quantum algorithms to solve DLP in finite groups, however, nothing is known to solve DLP in group rings or non-abelian groups $[3]$.  Therefore, we can say that DLP in  group rings is safer than DLP in groups, even in quantum word in the sense that no known literature is available.   Hence, we conclude that   ElGamal type group ring public key  cryptosystem is  more secure than the existing ElGamal  cryptosystem.
\begin{remark}
For ElGamal type group ring public key cryptosystem, $R$ can be ring of integers also. In that case too,  we  write the message in blocks of length $t$ with $|G|\geq t$. Example involving integral group ring is discussed in next section.
\end{remark}
\section{Examples}In this section, we discuss some examples of the group ring public key cryptosystems  defined in the preceding section. From now onwards,  we write  GRPKC as a short form for group ring public key cryptosystem.  Most of the results including determination of units, computation of powers of units and their corresponding inverses can be obtained  using  GAP (Groups, Algorithm and Programming) $[7]$, MAGMA or MATLAB.
\subsection{Elliptic ElGamal type GRPKC} For the better understanding of Elliptic ElGamal type GRPKC, we give a simple example by considering the small parameters. Suppose that the message $m$ of  arbitrary size consists of only English alphabets (this restriction is made only for understanding), let's say $\{$A, B, $\cdots$, Z$\}$. So, we map these alphabets on the elements of $\mathbb{F}_{29}$ by the mapping $$\{\text{A}\to 0, \ \text{B}\to 1, \ \text{C}\to 2, \ \cdots,\ \text{Z}\to 25\}.$$Suppose the message is ``ARMY". We consider $G=C_{29}=\ <g>$, i.e. cyclic group of order $29$. In terms of element of group ring, we write the message in the form $$r=0e+17g+12g^2+24g^3=\big[  0\ \ \ 17 \ \ \ 12 \ \ \ 24\big].$$Consider the elliptic curve $$E(\mathbb{F}_{29}):y^2=x^3+4x+20\mod 29$$ and  a point $P=(8, 10)$ on it. We consider the trivial unit $g$ of $\mathbb{F}_{29}G$ and choose the private key $(n_1, g^{-n_2})=(4, g^{-3})$. Public key is $(P, Q, A)$ where $$Q=(6, 17),\ \ A= g^3.$$ Let $n_3=3$ and  using Theorem $2.1$, we get $$n_3Q=(3, 28), \ \ 2n_3Q=(24, 22), \ \ 3n_3Q=(8, 19), \ \ 4n_3Q=(5, 22).$$Now we find the ciphertext $(C_1, C_2)$ corresponding to message $r$. Using $n_3=3$, we get $C_1=3P=(16,2)$. Further, we have $$(r\oplus n_3Q)= \big[ 
0\ \ \ 17\ \ \ 12\ \ \ 24\big]\oplus (3, 28)\hspace{50mm}$$ $$ = \big[0+3+28\ \ \ 17+24+22\ \ \  12+8+19 \ \ \ 24+5+22 \big] $$ $$=\big[ 2\ \ \ 5\ \ \ 10\ \ \ 22\big].$$ Above yields $$C_2=(r+n_3Q)*A= 2g^3+5g^4+10g^5+22g^6=\big[
0\ \ \ 0\ \ \ 0\ \ \ 2\ \ \ 5\ \ \ 10\ \ \ 22
\big].$$Therefore the ciphertext is $$(C_1, C_2)= \big((16,2), \big[
0\ \ \ 0\ \ \ 0\ \ \ 2\ \ \ 5\ \ \ 10\ \ \ 22
\big]\big).$$For decryption,  Alice first of all use the second part of her private key to get $$C_2*g^{-3}= \big[
0\ \ \ 0\ \ \ 0\ \ \ 2\ \ \ 5\ \ \ 10\ \ \ 22
\big]*g^{-3}= (2g^3+5g^4+10g^5+22g^6)*g^{-3}$$ $$= 2+5g+10g^2+22g^3=\big[
2\ \ \ 5\ \ \ 10\ \ \ 22\big].$$ Now Alice use  first part $n_1$ of her private key to obtain $$n_1C_1= 4(16,2)= (3, 28)$$ and  multiples of $(3, 28)$. Clearly Alice needs to compute only $4$ multiples of $(3, 28)$ as the remaining terms are zero in $C_2*g^{-3}$. Decrypted message is $$(C_2*g^{-3})\ominus (3, 28)=\big[
2\ \ \ 5\ \ \ 10\ \ \ 22\big]\ominus (3, 28) \hspace{50mm}$$ $$ = \big[
2-3-28\ \ \  5-24-22\ \ \  10-8-19 \ \ \ 22-5-22 \big] $$ $$=\big[0\ \ \ 17\ \ \ 12\ \ \ 24\big]= \text{ARMY}.$$
We have chosen trivial unit for this example so that it can be done manually. In the upcoming  examples, we consider non-trivial units which actually make the scheme worthy.
\subsection{Elliptic ElGamal type GRPKC for $R=\mathbb{Z}$ and $G=D_{10}$} 
In this subsection, we provide an example which highlights the importance of Remark $2.3$. We consider the integral group ring of group $G=D_{10}$ which means  message has integer values. $G$ can be represented as $$G=\langle r, s\ | \ r^5=1=s^2, sr=r^4s\rangle=\{1, r, r^2, r^3, r^4, s, rs, r^2s, r^3s, r^4s\}.$$  Consider the element $u=1-2r-2r^2+r^3+3r^4$ of $RG$. It can be verified that $u$ is a unit of $RG$ with  $u^{-1}=-2+3r-2r^2+r^3+r^4$. Now, consider the elliptic curve$$E(\mathbb{F}_{263}):y^2=x^3+2x+3\mod 263$$ and  a point $P=(200, 39)$ on it. Choose the private key $(10, u^{-8})$ with the corresponding public key $(P, Q, A)$ where $$Q=10P= (47, 78), \ A= u^8= -1576239+602070r+1948339r^2+602070r^3-1576239r^4$$ and $$u^{-8}= 602070+602070r-1576239r^2+1948339r^3-1576239r^4$$
Let $n_3= 5$ and using Theorem $2.1$, we get $$n_3Q=(180, 115), \ \ 2n_3Q=(5, 123), \ \ 3n_3Q=(128, 81), \ \ 4n_3Q=(127, 213)\ \ 5n_3Q=(17, 189)$$ $$6n_3Q=(102, 90), \ \ 7n_3Q=(74, 155), \ \ 8n_3Q=(142, 89), \ 9n_3Q=(144, 228)\ \ 10n_3Q=(139, 12).$$Now we find the ciphertext $(C_1, C_2)$ corresponding to message $$r'= \big[
1\ \ \ 2\ \ \ -1\ \ \ 6\ \ \ 0\ \ \ 8\ \ \ 0\ \ \ 3\ \ \ 9\ \ \ -5\big].$$ Using $n_3=5$, we get $C_1=5P=(251,155)$. Further, we have $$(r'\oplus n_3Q)= \big[
1\ \ \ 2\ \ \ -1\ \ \ 6\ \ \ 0\ \ \ 8\ \ \ 0\ \ \ 3\ \ \ 9\ \ \ -5\big] \oplus (180, 115)\hspace{37mm}$$ $$ = \big[
296\ \ \ 130\ \ \  208\ \ \  346\ \ \ 206\ \ \  200\ \ \ 229\ \ \ 234\ \ \ 381\ \ \  146\big]$$ Above yields $$C_2=(r'\oplus n_3Q)*A= 251904460+255117992r-94232542r^2-313356578r^3-99432146r^4$$ $$+277795332s+294897228rs-95538493r^2s-353942935r^3s-123209942r^4s .$$Therefore the ciphertext is $(C_1, C_2)$. Decryption can be done similarly as shown in subsection $3.1$. Clearly, this example shows the beautiful mixing of points of elliptic curve and the  message.
\newline
Now we discuss an example as said in Remark $2.4$. 
\subsection{ Elliptic ElGamal type GRPKC   for $R=\mathbb{F}_2$ and $G=D_{10}$}   We assume that message is in bits or in other words $R=\mathbb{F}_2$, $G$ is  Dihedral group of order $10$ with the same representation considered in last example.  Clearly, the group ring $RG$ has $1024$ elements. Consider the element $$u=r^2+r^4+s$$ of $RG$. It can be verified that $u$ is a unit of $RG$ with inverses $$u^{-1}=1+r+r^3+r^4+s+rs+r^4s.$$ Rest of the setting  including choice of elliptic curve, selection of point $P$ on elliptic curve, public key and  private key is exactly same to the one considered in last example. So, the private and public keys are  $(10, u^{-8})$ and  $(P, Q, A)$ respectively  where $$Q=10P= (47, 78), \ A= u^8= 1+r+r^2+rs+r^2s+r^3s+r^4s\ \ \text{and}\ \ u^{-8}= r^4+rs+r^4s.$$
Let $n_3= 5$ and using Theorem $2.1$ we can get the multiples of $n_3Q$ as done in last example.  Now we find the ciphertext $(C_1, C_2)$ corresponding to message $$r'= \big[
1\ \ \ 0\ \ \ 0\ \ \ 1\ \ \ 0\ \ \ 1\ \ \ 1\ \ \ 0\ \ \ 0\ \ \ 1\big].$$ Using $n_3=5$, we get $C_1=5P=(251,155)$. Further, we have $$(r'\oplus n_3Q)= \big[
1\ \ \ 0\ \ \ 0\ \ \ 1\ \ \ 0\ \ \ 1\ \ \ 1\ \ \ 0\ \ \ 0\ \ \ 1\big] \oplus (180, 115)\hspace{37mm}$$ $$ = \big[
296\ \ \ 128\ \ \  209\ \ \  340\ \ \ 206\ \ \  192\ \ \ 229\ \ \ 231\ \ \ 372\ \ \  151\big]\text{mod}\ 2$$ $$=  \big[
0\ \ \ 0\ \ \  1\ \ \  0\ \ \ 0\ \ \  0\ \ \ 1\ \ \ 1\ \ \ 0\ \ \  1\big].$$Above yields $$C_2=(r'\oplus n_3Q)*A= 1+r^2+r^4+s+rs+r^4s= \big[
1\ \ \ 0\ \ \  1\ \ \  0\ \ \ 1\ \ \  1\ \ \ 1\ \ \ 0\ \ \ 0\ \ \  1\big].$$Therefore the ciphertext is $(C_1, C_2)$. For decryption,  Alice first of all use the second part of her private key to get $C_2*A^{-1}$. Now Alice use  first part $n_1$ of her private key to obtain $n_1C_1=  (180, 115)$ and  its multiples. Clearly Alice needs to compute only $10$ multiples of $(180, 115)$ which is decided by the length of group. These multiples are then added to $C_2*A^{-1}$ under modulo $2$ to get the original message.

Now we discuss some examples for the feel of ElGamal GRPKC. All the three examples discussed above involve small parameters, but for the practical implementation of these schemes, we need to work with large numbers or in other words, large power of units. In the next example, we take care of this thing and consider  large (reasonably) parameters. 
\subsection{ ElGamal type GRPKC  for an abelian group}  Let $r=11110100001$ is the message, $R=\mathbb{F}_2$ and $G$ is cyclic group of order $11$, i.e.  $$G=\{g^i:0\leq 1\leq 10\}$$ and we denote $g_i=g^i$. Clearly, the group ring $RG$ has $2048$ elements. Given message as an element of group ring is $$r=1+g+g^2+g^3+g^5+g^{10}= \big[
1\ \ \ 1\ \ \ 1\ \ \ 1\ \ \  0\ \ \  1\ \ \  0\ \ \  0\ \ \  0\ \ \  0 \ \ \  1
\big].$$
 Take a  unit  $u=1+g+g^3$  (verify that $u^{-1}=1+g+g^3+g^6+g^7+g^8+g^{10}$) having order $1023$ and calculate $$A_1=u^{400}=g^2+g^4+g^6+g^7+g^{10}.$$   Further, choose another unit $$v=1+g+g^2+g^4+g^5+g^6+g^8+g^9+g^{10}\quad \text{with}\quad 
v^{-1}=1+g+g^3+g^6+g^9$$ and a secret positive integer $n_2=33$.  Public key is $(A_1, A_2,v)$ with $A_2=v^{33}.$ Private key is $(A_1^{-1}, 33)$ with $$A_1^{-1}=1+g+g^2+g^5+g^7+g^8+g^{10}.$$
Using  public key the obtained ciphertext with ephemeral key $k=19$ is $(C_1, C_2)$ where $$C_1=v^{19}=g^5+g^7+g^9$$ and $$C_2= (r*A_1)*A_2^k= (r*u^{400})*v^{627}=1+g^2+g^4+g^6+g^9+g^{10}$$ $$=\big[
1\ \ \ 0\ \ \ 1\ \ \ 0\ \ \  1\ \ \  0\ \ \  1\ \ \  0\ \ \  0\ \ \  1 \ \ \  1
\big].$$
After obtaining the ciphertext $(C_1, C_2)$, Alice computes $(C_2*C_1^{-n_2})*A_1^{-1}$  using her private key   to get the message.
\newline 
Next example is related to  Integral ElGamal GRPKC. 
\subsection{ ElGamal type GRPKC   for $R=\mathbb{Z}$ and $G=C_8$}  Let $RG$ be the integral group ring with $G=C_8=\langle x\rangle$, cyclic group of order $8$. We denote $g_i=x^i, 0\leq i\leq 7$ and consider the message blocks of length $8$ for encryption. Let $$r=\big[
0\ \ \ 1\ \ \ -10\ \ \ 17\ \ \ 0\ \ \ 0\ \ \ 0\ \ \ -4\big]$$ is the message to be encrypted. In terms of element of group ring, $r$ can be written as $$r=x-10x^2+17x^3-4x^7.$$ Here we choose small parameters otherwise the expressions becomes so large to print here  but one can choose arbitrary large numbers. Consider the element  $u=2e+x-x^3-x^4-x^5+x^7$ of $RG$. It can be verified that $u$ is a unit of $RG$ with   $u^{-1}=2e-x+x^3-x^4+x^5-x^7$ having infinite order. Take $n_1=2$ and calculate $$A_1=u^2=9e+6x-6x^3-8x^4-6x^5+6x^7 \quad 
\text{with}\quad A_1^{-1}= 9e-6x+6x^3-8x^4+6x^5-6x^7.$$ Further, consider another unit $$v=50e+35x-35x^3-49x^4-35x^5+35x^7\ \ \text{with}\ \ v^{-1}=50e-35x+35x^3-49x^4+35x^5-35x^7$$ of $RG$ and  secret positive integer $n_2=2$.  Public key is $(A_1, A_2, v)$ where $$A_2=v^2=9801e+6930x-6930x^3-9800x^4-6930x^5+6930x^7$$and  private key is $(A_1^{-1}, 2)$.
Using  public key the obtained ciphertext with ephemeral key $k=1$ is $(C_1, C_2)$ where $C_1=v$ and $$C_2= (r*A_1)*v^2=-4708320e-2021231x+1849862x^2+4637345x^3+4708320x^4$$ $$+2021232x^5-1849872x^6-4637332x^7.$$
After obtaining the ciphertext $(C_1, C_2)$, Alice computes $(C_2*v^{-2})*u^{-1}$  using her private key   to obtain the message. 

We end this section by giving an example in which we set a public key which is multiplication of Bass cyclic unit and Bicyclic unit of integral group ring. This public key or unit has an advantage that it is not of the known type and therefore it becomes computationally infeasible to obtain the inverse of such a unit.
\subsection{ ElGamal type GRPKC  for $R=\mathbb{Z}$ and $G=S_5$} 

 Let $RG$ be the integral group ring with $G=S_5$, i.e. symmetric group of degree $5$. We denote $G=\{g_i, 0\leq i\leq 119\}$ (for any order of our choice)  and consider the message blocks of length $120$ for encryption. Now before the setting of public key,  we recall the definitions of Bass cyclic and Bicyclic units. \begin{definition}
Bass cyclic unit: Let $g$ be an element of a group $G$ having order $n$. Choose an integer $i$ coprime to $n$ such that $1<i<n-1$. Then the element $$u=(1+g+g^2+\cdots+g^{i-1})^{\phi(n)}+\frac{1-i^{\phi(n)}}{n}\hat{g}$$ is a unit of integral group ring $RG$. Here $\phi(n)$ denotes the Euler's totient function and $$\hat{g}=1+g+g^2+\cdots+g^{n-1}.$$
\end{definition}

\begin{definition}
Bicyclic unit: Let $g, h$ be two elements of a group $G$ such that $|g|=n$.  Then the element $$u=1+(g-1)h\hat{g}$$ is a bicyclic unit of integral group ring $RG$ corresponding to $g$ and $h$ where $\hat{g}$ has the same meaning as in Definition $3.1$.
\end{definition} For additional reading on these units see $[15]$. Let $g=(1\ 4\ 2\ 5\ 3)$ be an element of $S_5$ having order $5$ and $i=3$. Then the  Bass cyclic unit $x$ corresponding to $g$ and $i$ is \small{$$x=1*(1)-2*(1\ 2\ 3\ 4\ 5)+3*(1\ 3\ 5\ 2\ 4)-2*(1\ 4\ 2\ 5\ 3)+1*(1\ 5\ 4\ 3\ 2).$$}
\normalsize{Further, consider elements $(1\ 2\ 3\ 4\ 5)$ and $(1\ 2)$ of $S_5$ and use these to obtain the corresponding bicyclic unit} 
\small{ $$y=1*(1)+1*(2\ 3\ 4\ 5)-1*(2\ 5\ 4\ 3)-1*(1\ 2)+1*(1\ 2\ 4)(3\ 5)-1*(1\ 3\ 4\ 5)+1*(1\ 3)(2\ 5\ 4)$$ $$+1*(1\ 4\ 3\ 2)-1*(1\ 4)(2\ 3\ 5)+1*(1\ 5)-1*(1\ 5\ 3)(2\ 4).\hspace{50mm}$$}\normalsize{Let $u=x*y$ where} \small{$$\hspace{3mm} u=1*(1)-5*(4\ 5)+5*(3\ 4)-3*(2\ 3)+3*(2\ 3\ 4\ 5)-5*(1\ 2\ 3\ 4)-2*(1\ 2\ 3\ 4 \ 5)+5*(1\ 2\ 3\ 5)$$ $$\hspace{5mm}-3*(1\ 2\ 4\ 5)+3*(1\ 2\ 4)(3\ 5)+5*(1\ 3)(2\ 4\ 5)+3*(1\ 3\ 5\ 2\ 4)-5*(1\ 3\ 5)(2\ 4)+3*(1\ 3)(2\ 5\ 4)$$ $$\hspace{5mm}-3*(1\ 3\ 4)(2\ 5)+3*(1\ 4 \ 3 \ 2)-3*(1\ 4\ 2)(3\ 5)-2* (1\ 4 \ 2\ 5 \ 3)-5*(1\ 4\ 3)(2\ 5)+5*(1\ 4)(2\ 5\ 3)$$ $$+1*(1\ 5\ 4\ 3\ 2)-5*(1\ 5\ 3\ 2)+5*(1\ 5\ 4\ 2) -3*(1\ 5\ 4\ 3)+3*(1\ 5)\hspace{40mm}$$} \normalsize{with} 
\small{$$\hspace{3mm} u^{-1}=-2*(1)-1*(2\ 3\ 4\ 5)+1*(2\ 5\ 4\ 3)+1*(1\ 2)-2*(1\ 2\ 3\ 4\ 5)-1*(1\ 2\ 4)(3 \ 5)+1*(1\ 3\ 4\ 5)-1*(1\ 3)( 2\ 5\ 4)$$ $$\hspace{10mm}+1*(1\ 3\ 5\  2\ 4)-1*(1\ 4\ 3\ 2) +1*(1\ 4)(2\ 3\ 5)+3*(1\ 4\ 2\ 5\ 3)+1*(1\ 5\ 4\ 3\ 2)-1*(1\  5)+1*(1\ 5 \ 3) (2, 4).$$}\normalsize{Choose another unit} \small{ $$v=-15*e+6*(1\ 2\ 3\ 4\ 5)+6*(1\ 3\ 5\ 2\ 4)-15*(1\ 4\ 2\ 5\ 3)+19*(1\ 5\ 4\ 3\ 2)$$} \normalsize{with} \small{ $$\hspace{3mm}v^{-1}=6*e+19*(1\ 2\ 3\ 4\ 5)+6*(1\ 3\ 5\ 2\ 4)-15*(1\ 4\ 2\ 5\ 3)-15*(1\ 5\ 4\ 3\ 2).$$}
\normalsize{ Let the message to be encrypted is}\small{ $$r=-1*(2\ 3\ 4\ 5)+2*(1\ 2\ 4)(3\ 5)-3*(1\ 3)(2\ 5\ 4)+4*(1\ 4\ 3\ 2)-5*(1\ 5).$$}
 \normalsize{Take $n_1=1$ and $n_2=3$ which means public key is $(A_1=u, A_2=v^3, v)$ where} \small{ $$v^3=12816*e-33552*(1\ 2\ 3\ 4\ 5)+41473*(1\ 3\ 5\ 2\ 4)-33552*(1\ 4\ 2\ 5\ 3)+12816*(1\ 5\ 4\ 3\ 2).$$}\normalsize{Corresponding  private key is $(A_1^{-1}, 3)$. 
Using  public key the obtained ciphertext with ephemeral key $k=2$ is $(C_1, C_2)$ where} \small{ $$C_1=273*e+273*(1\ 2\ 3\ 4\ 5)-714*(1\ 3\ 5\ 2\ 4)+883*(1\ 4\ 2\ 5\ 3)-714*(1\ 5\ 4\ 3\ 2)$$}\normalsize{ and}\small{ $$ C_2= (r*A_1)*(v^3)^2= 
-56*(1)+40*(3\ 4\ 5)-228001425142*(2\ 3\ 4\ 5)-26*(2\ 3\   5)-9*( 2\ 4\ 5)+51*( 2\  4)(3 \ 5)$$ $$\hspace{8mm}-56*(1\ 2\ 3\ 4\ 5)+40*(1\ 2\ 3\ 5\ 4)-26*(1\ 2\ 4 \ 5 \ 3)+51*(1\ 2\ 5\ 4\ 3)-9*(1\ 2\ 5\ 3\ 4)+51166299748*(1\ 2\ 4)$$ $$\hspace{8mm}(3\ 5)+51*(1\ 3\ 2)-9*(1\ 3\ 5\ 4\ 2)+40*(1\ 3)( 2\ 4)-56*(1\ 3\ 5\ 2\ 4)+145212613069*(1\ 3)(2\ 5\ 4)-26*(1\ 3\ 2\ $$ $$\hspace{10mm}5\ 4)-26*(1\ 4\ 2)-286125243291*(1\ 4 \ 3 \ 2)-9*(1\ 4\ 3)+51* (1\ 4 \ 5)-56*(1\ 4\ 2\ 5\ 3)+40*(1\ 4\ 3 \ 2\ 5)-$$ $$\hspace{10mm}56*(1\ 5\ 4\  3\ 2)+40*(1\ 5\ 2)+317747755613*(1\ 5)-26*(1\ 5)(3\ 4) -9*(1\ 5)(2\ 3)+51*(1\ 5\ 2\ 3\ 4)
.$$}
\normalsize{After obtaining the ciphertext $(C_1, C_2)$, Alice  can easily obtain the message.}
\section{Units of group rings using GAP}
 Both the encryption schemes of proposed cryptosystems involve units of group rings and therefore in this section we discuss about the units of group rings. 
  There are only few recipes for the construction of units but still   a vast literature is available on units of group rings.  Using LAGUNA package  of GAP $[5]$, group $V(RG)$ of  normalized units of group algebra $RG$ with $R=\mathbb{F}_p$ for some prime $p>0$ and  a finite $p$-group can be computed efficiently for smaller sizes of $RG$. Further,  as the unit group $U(RG)$ of $RG$ is $\mathbb{F}_p^*\times V(RG)$ which means we can efficiently compute the complete unit group. But with the increase in size of $|G|$, GAP is very inefficient, for instance, the computation of number of elements and generators of normalized unit group of $V(\mathbb{F}_{5}G)$ where $G$ is a cyclic group of order $5^5$ takes $38$ minutes (approx)   and program may crash for groups of size greater than $5000$. Despite of this drawback,  GAP still covers a  large class of unit groups.  
 For finite groups other than $p$-groups, units can be obtained using the Wedderga package of GAP $[2]$ with the additional restrictions that $R$ is semisimple  and $|G|$ is invertible in $R$. Basically Wedderga package is not for the calculation of units directly. Its main purpose is to find the Wedderburn decomposition $[15]$ of a group algebra. From this wedderburn decomposition, we can easily obtain the structure of unit groups.  Without the use of LAGUNA and Wedderga we  can also generate the units of integral group rings as well as any modular group ring  using GAP.
 \section{Conclusion}
From security point of view, we already discussed that both the proposed cryprosystems have greater security than that of existing cryptosystems. In the coming era of Quantum computers, this enhanced security would play a great role because no such quantum algorithm is known to solve DLP in group rings.  Novelty of this work is the involvement of group rings and their  units  in the field of public key cryptography.  As mentioned there is a huge literature available on the structure of unit group of group rings but very few articles are there in which unit group is computed explicitly in terms of elements of group ring. Since the latter is more important for frequent use of group rings in cryptography, much more efforts are required for the explicit representation of the unit group in terms of the elements of  group ring. 
\bibliographystyle{plain}

\begin{thebibliography}{10}
\small{

\bibitem {}{ I. Anshel, M. Anshel, D. Goldfeld,} \newblock{\em An algebraic method for public-key
cryptography,}  \newblock{Math. Res. Lett., $6$, $287$-$291$, 1999.}



\bibitem{}{G. K. Bakshi, O. B. Cristo, A. Herman, A. Konovalov, S. Maheshwary, A. Olteanu, A. Olivieri, \'A. del R\'io and I. V. Geldar}, \newblock{\em Wedderga - wedderburn decomposition of group algebras-a GAP package}, \newblock{  https://gap-packages.github.io/wedderga/, Version $4.9.5$, $2018.$}

\bibitem{}{D. J. Bernstein, J. Buchmann, E. Dahmen}, \newblock{\em Post quantum cryptography}, \newblock{Springer, $2009$.}

\bibitem{}{I. F. Blake, G. Seroussi, N. P. Smart,} \newblock{\em Elliptic curves in cryptography}, \newblock{volume $265$ of London mathematical society lecture note series, Cambridge university press, $2000$.}


\bibitem{}{V. Bovdi, A. Konovalov, R. Rossmanith, 
C. Schneider}, \newblock{\em LAGUNA - Lie algebras and units of group Algebras}, \newblock{ https://gap-packages.github.io/laguna/, Version $3.8.0$, $2017.$}

\bibitem{}{M. Chandrashekhar,  L. Xiangxue}, \newblock{\em New efficient key authentication protocol for public key cryptosystem using DL over multiplicative group}, \newblock{Journal of Information and Optimization Sciences, 39(2), 391-400, $2018$.}

\bibitem{}{GAP}, \newblock{\em  Groups, Algorithms, Programming}, \newblock{https://www.gap-system.org.}


\bibitem{} \newblock{ D. Garber}, \newblock{\em  Braid group cryptography}, https://arxiv.org/abs/0711.3941.


\bibitem{}{   J. Hoffstein, J. Pipher, J. H. Silverman,} \newblock{\em An introduction of mathematical cryptography}, \newblock{Springer, $2008$.}



\bibitem{}{  B. Hurley, T. Hurley,} \newblock{\em Group ring cryptography}, \newblock{Int. J. Pure Appl. Math., $69(1)$,  67-86, 2011.}


\bibitem{}{T. Hurley,} \newblock{\em Group rings for communications}, \newblock{Int. J. Group Theory, $4(4)$,  1-23, 2015.}


\bibitem{}{ N. Koblitz,} \newblock{\em Introduction to elliptic curves and modular forms}, \newblock{Springer, $1993$.}

\bibitem{}{   N. Koblitz,} \newblock{\em A course in number theory and cryptography}, \newblock{Springer, $1994$.}




\bibitem{}{A. Menezes,} \newblock{\em Elliptic curves public key cryptosysytems}, \newblock{The Kluwer international series in engineering and computer science, $234$, Kluwer academic publishers, $1993$.}

\bibitem{}{ C. P. Milies, S. K. Sehgal,} \newblock{\em An Introduction to group rings}, \newblock{Kluwer, $(2002)$.}


\bibitem{}{J. Proos, C. Zalka}, \newblock{\em Extension of Shor's to ECDSA}, \newblock{https://arxiv.org/pdf/quant-ph/0301141.pdf, $2003$.}


\bibitem{}{C. Shannon,} \newblock{\em Communication theory of secrecy systems}, \newblock{Bell systems technical journal, $4$, $1949$.}





\bibitem{}{P. W. Shor}, \newblock{\em Algorithms for quantum computation: discrete logarithms and factoring}, \newblock{Proceedings 35th Annual Symposium on Foundations of Computer Science,  IEEE Comput. Soc. Press., $1994$.}


\bibitem{}{  L. Shpilrain, G. Zapata,} \newblock{\em  Combinatorial group theory and public key cryptography},\newblock{ Appl. Algebra Engrg. Comm. Comput., $17(3$-$4)$,  291-302, 2006.}

\bibitem{}{W. Stallings,} \newblock{\em Cryptography and network security}, \newblock{Pearson, $2017$.}



\bibitem{}{P. Sundarayya and G. V. Prasad}, \newblock{\em A public key cryptosystem using Affine Hill Cipher under modulation of prime number}, \newblock{Journal of Information and Optimization Sciences, 40(4), 919-930, $2019$.}



\small}


 \end{thebibliography}

\end{document}